\documentclass[12pt,a4paper]{article}
\usepackage[latin2]{inputenc}
\usepackage{amsmath}
\usepackage{amsfonts}
\usepackage{amssymb}
\usepackage{graphicx}
\usepackage[left=2.00cm, right=2.00cm, top=2.50cm, bottom=2.50cm]{geometry}
\usepackage[T1]{fontenc}

\def\Q{{\mathbb Q}}
\def\Z{{\mathbb Z}}

\def\O{{\cal O}}

\newtheorem{lemma}{Lemma}
\newtheorem{theorem}[lemma]{Theorem}

\title{
Monogenity in totally complex sextic fields, revisited
}
\author{
Istv\'{a}n Ga\'{a}l \\
{\small University of Debrecen, Mathematical Institute} \\
{\small H--4002 Debrecen Pf.400., Hungary,} \\
{\small e--mail: gaal.istvan@unideb.hu},
}

\begin{document}
\baselineskip=17pt

\maketitle
\thispagestyle{empty}

\renewcommand{\thefootnote}{\arabic{footnote}}
\setcounter{footnote}{0}


\begin{abstract}
In addition to rather complicated general methods it is interesting and valuable to develop fast 
efficient methods for calculating generators of power integral bases in special types of number 
fields. We consider sextic fields containing a real cubic and a complex quadratic fields.
We develop a very simple and very efficient method to calculate generators of power integral
bases in this type of fields. Our method can be applied to infinite families of number fields, as
well. We substantially improve the former methods.
Our algorithm is illustrated with detailed examples, involving infinite 
parametric families.
\end{abstract}

\newpage

\section{Introduction}

In the following we shall denote 
by $\Z_K$ and $D_K$ the ring of integers and the discriminant, respectively,
of any number field $K$.

There is an extensive literature of monogenity of number fields and power integral bases,
see \cite{nark}, \cite{book}. A number field $K$ of degree $n$ is {\it monogenic} if $\Z_K$
is a simple ring extension of $\Z$, that is there exists $\alpha\in\Z_K$ with 
$\Z_K=\Z[\alpha]$. In this case $(1,\alpha,\ldots,\alpha^{n-1})$ is an integral basis of $K$,
called {\it power integral basis}. (We also call $\alpha$ the {\it generator}
of this power integral basis.) 
The algebraic integer $\alpha$ generates a power integral basis
if and only if its {\it index}
\[
I(\alpha)=\sqrt{\left|\frac{D(\alpha)}{D_K}\right|}
\]
is equal to 1, where $D(\alpha)$ is the discriminant of $\alpha$.

The calculation of generators of power integral bases can be reduced to certain
diophantine equations, called {\it index form equations}, cf. \cite{book}.

There exist general algorithms for solving index form equations in cubic, quartic, 
quintic, sextic fields, however the general algorithms for quintic and sextic fields 
are already quite tedious, see \cite{s6}.  
Therefore it is worthy to develop efficient methods for
the resolution of special types of higher degree number fields.

In this paper we study totally complex sextic fields that are composites of 
a totally real cubic and an imaginary quadratic fields.
These fields were investigated in \cite{totcomplexsextic} where we reduced the
relative Thue equation involved to absolute Thue inequalities.
That method was further developed in \cite{realthueimquadr}.
However some ideas of \cite{inhomrelthue} lead to a considerable improvement 
of that algorithm, what we are going to detail here.
We also note that a parametric family of this type of number fields, 
consisting of composites of  the simplest cubic fields and imaginary quadratic fields
was studied in \cite{gr}, but applying results on the connected simplest family
of relative Thue equations, which counts as a more complicated approach.

\section{Composites of real cubic and imaginary quadratic fields}

Let $\vartheta=\vartheta^{(1)},\vartheta^{(2)},\vartheta^{(3)}$ be the roots 
of the totally real polynomial $f(x)=x^3+a_2x^2+a_1x+a_0\in\Z[x]$
and let $L=\Q(\vartheta)$.
Let $0<d\in\Z$ be a square-free integer, set $M=\Q(i\sqrt{d})$.
Our purpose is to calculate all generators of power integral bases in
$K=L\cdot M=\Q(\vartheta,i\sqrt{d})$. Set
\[
\omega=i\sqrt{d}\; {\rm if}\; -d\equiv 2,3\; (\bmod{4})
\;\; {\rm and}\;\; \omega=\frac{1+i\sqrt{d}}{2}\; {\rm if} -d\equiv 1\;  (\bmod{4}).
\]
Denote by $\gamma'$ the conjugate of of any $\gamma\in M$.

To make our presentation as simple as possible, we formulate our statements for the order 
\[
\O=\Z[1,\vartheta,\vartheta^2,\omega,\omega\vartheta,\omega\vartheta^2].
\]
The cubic field $L$ very often happens to have integral basis $(1,\vartheta,\vartheta^2)$,
and if $D_L$ is relatively prim to $D_M$,
then indeed $\O=\Z_K$. However, otherwise our statements are applicable with minor modifications, 
see Remark 2.

Let us represent any $\alpha\in\O$ in the form
\begin{equation}
\alpha=x_0+x_1\vartheta+x_2\vartheta^2+y_0\omega+y_1\omega\vartheta+y_2\omega\vartheta^2=
X_0+X_1\vartheta+X_2\vartheta^2,
\label{alpha}
\end{equation}
where $x_j,y_j\in\Z,\; X_j=x_j+\omega y_j\in\Z_M\; (j=0,1,2)$.

The conjugates of $\alpha$ are obtained the following way:
\[
\alpha^{(j,1)}=x_0+x_1\vartheta^{(j)}+x_2(\vartheta^{(j)})^2+y_0\omega
+y_1\omega\vartheta^{(j)}+y_2\omega(\vartheta^{(j)})^2=
X_0+X_1\vartheta^{(j)}+X_2(\vartheta^{(j)})^2,
\]
\[
\alpha^{(j,2)}=x_0+x_1\vartheta^{(j)}+x_2(\vartheta^{(j)})^2+y_0\omega'
+y_1\omega'\vartheta^{(j)}+y_2\omega'(\vartheta^{(j)})^2=
X_0'+X_1'\vartheta^{(j)}+X_2'(\vartheta^{(j)})^2, 
\]
for $j=1,2,3$.

Keeping the coefficients as variables, consider the symmetric polynomial
\[
F(x_1,x_2,y_0,y_1,y_2)=
\]
\[
=(\alpha^{(1,1)}-\alpha^{(2,2)})(\alpha^{(1,1)}-\alpha^{(2,3)})
(\alpha^{(1,2)}-\alpha^{(2,1)})(\alpha^{(1,2)}-\alpha^{(2,3)})
(\alpha^{(1,3)}-\alpha^{(2,1)})(\alpha^{(1,3)}-\alpha^{(2,2)}),
\]
having all coefficients in $\Z$.

\begin{theorem}
\label{th1}
If $\alpha\in\O$ (in the representation (\ref{alpha})) generates a power integral basis
in $\O$, then the coefficients $x_1,x_2,y_0,y_1,y_2\in\Z$ of $\alpha$ satisty
\begin{equation}
F(x_1,x_2,y_0,y_1,y_2)=\pm 1,
\label{F3}
\end{equation}
and 
\begin{equation}
N_{L/Q}(y_0+y_1\vartheta+y_2\vartheta^2)=\pm 1.
\label{L}
\end{equation}
Further, if $-d\equiv 2,3\; (\bmod{4})$, then
\begin{equation}
|N_{L/Q}(x_1-(a_2+\vartheta)x_2)|\leq 1
\label{x12}
\end{equation}
and
\begin{equation}
|N_{L/Q}(y_1-(a_2+\vartheta)y_2)|\leq \frac{1}{(\sqrt{d})^3}.
\label{y12}
\end{equation}
If $-d\equiv 1\; (\bmod{4})$, then
\begin{equation}
|N_{L/Q}((2x_1+y_1)-(a_2+\vartheta)(2x_2+y_2)|\leq 8
\label{x128}
\end{equation}
and
\begin{equation}
|N_{L/Q}(y_1-(a_2+\vartheta)y_2)|\leq \frac{8}{(\sqrt{d})^3}.
\label{y128}
\end{equation}
\end{theorem}

\vspace{1cm}

Recall, that $a_2$ is the coefficient of the quadratic term of the defining polynomial of
$\vartheta$.

\vspace{1cm}

\noindent
{\bf Proof of Theorem \ref{th1}}.\\
The discriminant of the basis $(1,\vartheta,\vartheta^2,\omega,\omega\vartheta,\omega\vartheta^2)$
of $\O$ is
\begin{equation}
D_{\O}=D(\vartheta)^2\cdot D_M^3.
\label{DO}
\end{equation}
For $j,k\in\{1,2,3\},j\ne k$ we have
\[
\alpha^{(1,j)}-\alpha^{(1,k)}=(\vartheta^{(j)}-\vartheta^{(k)})
(X_1+(\vartheta^{(j)}+\vartheta^{(k)})X_2)=
(\vartheta^{(j)}-\vartheta^{(k)})(X_1-(a_2+\vartheta^{(\ell)})X_2),
\]
\[
\alpha^{(2,j)}-\alpha^{(2,k)}=(\vartheta^{(j)}-\vartheta^{(k)})
(X_1'+(\vartheta^{(j)}+\vartheta^{(k)})X_2')=
(\vartheta^{(j)}-\vartheta^{(k)})(X_1'-(a_2+\vartheta^{(\ell)})X_2'),
\]
where $\ell=\{1,2,3\}\setminus \{j,k\}$.
Therefore
\begin{equation}
\prod_{s=1}^2\prod_{1\leq j<k\leq 3}(\alpha^{(s,j)}-\alpha^{(s,k)})
=D(\vartheta)\cdot N_{M/Q}(N_{K/M}(X_1-(a_2+\vartheta)X_2)).
\label{F1}
\end{equation}
Further,
\[
\alpha^{(1,j)}-\alpha^{(2,j)}=(\omega-\omega')(y_0+y_1\vartheta+y_2\vartheta^2),
\]
hence
\begin{equation}
\prod_{j=1}^3(\alpha^{(1,j)}-\alpha^{(2,j)})=
(\omega-\omega')^3 N_{L/Q}(y_0+y_1\vartheta+y_2\vartheta^2).
\label{F2}
\end{equation}
The remaining factors of $I(\alpha)$ are just those of $F(x_1,x_2,y_0,y_1,y_2)$.
In view of (\ref{DO}) this implies, that 
\[
I(\alpha)=\sqrt{\left|\frac{D(\alpha)}{D_{\O}}\right|}=
N_{M/Q}(N_{K/M}(X_1-(a_2+\vartheta)X_2))\cdot N_{L/Q}(y_0+y_1\vartheta+y_2\vartheta^2)
\cdot F(x_1,x_2,y_0,y_1,y_2),
\]
or equivalently, $I(\alpha)=1$, if and only if (\ref{F3}), (\ref{L}) and
\begin{equation}
N_{M/Q}(N_{K/M}(X_1-(a_2+\vartheta)X_2))=\pm 1,
\label{nn}
\end{equation}
simultaneously hold.

Our present improvement concerns this last equation. 
(\ref{nn}) implies 
\begin{equation}
|N_{K/M}(X_1-(a_2+\vartheta)X_2))|=\pm 1
\label{nna}
\end{equation}
since the norm $M/Q$ is just the product of the above norm and its complex conjugate.
We obtain 
\[
\left|\prod_{j=1}^3((x_1+\omega y_1)-(a_2+\vartheta^{(j)})(x_2+\omega y_2))\right|= 1.
\]
Set $\beta^{(j)}=(x_1+\omega y_1)-(a_2+\vartheta^{(j)})(x_2+\omega y_2),j=1,2,3$.
Obviously $|{\rm Re}(\beta^{(j)})|\leq |\beta^{(j)}|$ and 
$|{\rm Im}(\beta^{(j)})|\leq |\beta^{(j)}|$,
for $j=1,2,3$, whence
\begin{equation}
\prod_{j=1}^3|{\rm Re}(\beta^{(j)})|\le \prod_{j=1}^3|\beta^{(j)}|=1
\;\; {\rm and}\;\; 
\prod_{j=1}^3|{\rm Im}(\beta^{(j)})|\le \prod_{j=1}^3|\beta^{(j)}|=1.
\label{reim}
\end{equation}

If $-d\equiv 2,3\; (\bmod{4})$ then 
\[
{\rm Re}(\beta^{(j)})=x_1-(a_2+\vartheta^{(j)})x_2,\;\;
{\rm Im}(\beta^{(j)})=y_1-(a_2+\vartheta^{(j)})y_2,
\]
and (\ref{reim}) implies (\ref{x12}), (\ref{y12}).

If $-d\equiv 1\; (\bmod{4})$ then 
\[
{\rm Re}(\beta^{(j)})=\frac{1}{2}((2x_1+y_1)-(a_2+\vartheta^{(j)})(2x_2+y_2)),\;\;
{\rm Im}(\beta^{(j)})=\frac{\sqrt{d}}{2}(y_1-(a_2+\vartheta^{(j)})y_2),
\]
and (\ref{reim}) implies (\ref{x128}), (\ref{y128}).
\hfill $\Box$

\vspace{1cm}

\noindent
{\bf Remark 1.}\\
We already had equations (\ref{L}) and (\ref{nna}) in \cite{totcomplexsextic},
but (\ref{x12}), (\ref{y12}), (\ref{x128}), (\ref{y128}) are
much stronger than the corresponding inequalities of Theorem 2.2 in \cite{totcomplexsextic}.
These will be very useful in applications, see Sections \ref{appl}, \ref{ex}.

\vspace{1cm}

\noindent
{\bf Remark 2.}\\
If $\O$ is not equal to $\Z_K$, then any $\alpha\in\Z_K$ can be written in the form
\[
\alpha=\frac{x_0+x_1\vartheta+x_2\vartheta^2+y_0\omega+y_1\omega\vartheta+y_2\omega\vartheta^2}{g}
\]
with a common denominator $g\in\Z$. This implies
\[
N_{M/Q}(N_{K/M}(X_1-(a_2+\vartheta)X_2))\cdot N_{L/Q}(y_0+y_1\vartheta+y_2\vartheta^2)
\cdot F(x_1,x_2,y_0,y_1,y_2)=\pm g^{15}.
\]
This factor $g^{15}$ splits into a product of integers $g_1,g_2,g_3$ with
$g_1g_2g_3=g^{15}$ such that $g_1,g_2,g_3$ occur on the right hand sides of 
the above three factors,
that is we get 
(\ref{F3}) with right hand side $\pm g_1$,
(\ref{L}) with right hand side $\pm g_2$
and
(\ref{x12}), (\ref{y12}), (\ref{x128}), (\ref{y128}) with right hand sides 
equal to $g_1$ times the original right hand sides.

\vspace{1cm}

\section{How to apply Theorem \ref{th1}?}
\label{appl}

If $d\ne 1$ in case $-d\equiv 2,3\; (\bmod{4})$, then (\ref{y12}) has only the trivial solutions 
$y_1=y_2=0$, whence by (\ref{L}) we obtain $y_0=\pm 1$. 

Similarly, if $d\ne 3$ in case $-d\equiv 1\; (\bmod{4})$, then (\ref{y128}) has only the trivial 
solutions $y_1=y_2=0$, whence by (\ref{L}) we obtain $y_0=\pm 1$. 

For the Gaussian integers ($d=1$) and Euler integers ($d=3$), (\ref{y12}), resp. (\ref{y128})
are Thue inequalities with some small right hand sides.

For any given $y_1,y_2$ equation (\ref{L}) is just a cubic polynomial equation in the
integer variable $y_0$.

In case $-d\equiv 2,3\; (\bmod{4})$ we can determine $x_1,x_2$ from the Thue equation (\ref{x12}).

In case $-d\equiv 1\; (\bmod{4})$, given $y_1,y_2$ 
we can determine $x_1,x_2$ from the Thue equation (\ref{x128}),
solving it with right hand sides $0,1,\ldots,8$.

Remark that using Magma \cite{magma} or Kash \cite{kant} it is no problem to solve
cubic Thue equations with small right hand sides within a few seconds.

Having calculated $x_1,x_2,y_0,y_1,y_2$ we have to check if $\alpha$ of (\ref{alpha}) 
has indeed index 1 (equations (\ref{F3}), (\ref{L}), (\ref{nn}) together are equivalent with
$I(\alpha)=1$, but the inequalities (\ref{x12}), (\ref{y12}) and  (\ref{x128}), (\ref{y128}),
respectively, are weaker than (\ref{nn}).

It is easy to explicitly calculate the polynomial $F(x_1,x_2,y_0,y_1,y_2)$.
Equation (\ref{F3}) is very useful if we consider monogenity in infinite parametric
families of number fields, see Section \ref{ex}.

\section{Examples}
\label{ex}

\subsection{Example 1}

Let $t\in\Z$ and consider the infinite parametric family of fields
\begin{equation}
f_t(x)=x^3-(t^4-t)x^2+(t^5-2t^2)x+1.
\label{ft}
\end{equation}
According to \cite{bg} the polynomial $f_t$ has three real roots for $t\ge 2$.
In the following let $t\geq 2$ and denote by $\vartheta_t$ a root of $f_t$.
Let $L_t=\Q(\vartheta_t)$.

Set $F_t(x,y)=y^3f_t(x/y)$. The infinite parametric family of Thue equations
\begin{equation}
F_t(x,y)=N_{L_t/Q}(x-\vartheta_t y)=1\;\; {\rm in}\;\; x,y\in \Z
\label{Ft}
\end{equation}
was recently considered M. A. Bennett and A. Ghadermarzi \cite{bg}.
They showed that for $t\ne -1$ all solutions of the above equation are
$(x,y)=(1,0),(0,1),(t,1),(t^4-2t,1),(1-t^3,t^8-3t^5+3t^2)$ and for $t=-1$
it has the additional solution $(x,y)=(6,-5)$.

Let $d>1$ be a square-free integer with $-d\equiv 2,3\; (\bmod{4})$ and let 
$\omega=i\sqrt{d}$.

Consider the order 
$\O_{t,d}=\Z[1,\vartheta_t,\vartheta_t^2,\omega,\omega\vartheta_t,\omega\vartheta_t^2]$
of the number field $K_{t,d}=\Q(\vartheta_t,i\sqrt{d})$.
Remark that $(1,\vartheta_t,\vartheta_t^2)$ very often happens to be an integer basis 
of $\Q(\vartheta_t)$ (in case the discriminant of $f_t$ is square-free, but also in other cases).
Further if $(1,\vartheta_t,\vartheta_t^2)$ is an integer basis of $\Q(\vartheta_t)$ and
the discriminant of $f_t$ is co-prime to $4d$, then $\O_{t,d}$ 
is just the  ring of integers of $K_{t,d}$.

{\vspace{1cm}

We have

\begin{theorem}
\label{th2}
For $t\ge 2,d>1,-d\equiv 2,3\; (\bmod{4})$, the order $\O_{t,d}$ is never monogenic. 
\end{theorem}

{\vspace{1cm}

\noindent
{\bf Proof of Theorem \ref{th2}}.\\

\noindent
For $d>1$ we have $y_1=y_2=0$ from (\ref{y12}) and $y_0=\pm 1$ from (\ref{L}).
The solutions of (\ref{x12}) we obtain form the result of \cite{bg} on equation
(\ref{Ft}), by a suitable transformation (if $x,y$ is a solution of (\ref{Ft}),
then $x_1=x+a_2y,x_2=y$ is a solution of (\ref{x12}), where $a_2$ is the coefficient
of $x^2$ in the defining polynomial of $\vartheta_t$, that is, $a_2=t^4-t$ in our case).

We substitute all possible $x_1,x_2,y_0,y_1,y_2$ into $F(x_1,x_2,y_0,y_1,y_2)$. 
We obtain cubic polynomials in $d$. Their coefficients are polynomials in $t$.
All these coefficients are negative for $t\ge 2$, that is (\ref{F3}) can not be satisfied.
\hfill $\Box$

\subsection{Example 2}

Let $f(x)=x^3-x^2-2x+1$. This polynomial (with discriminant  49) has real roots.
Denote by $\vartheta$ a root of $f$. Let $L=\Q(\vartheta)$, with
integral basis $(1,\vartheta,\vartheta^2)$.

Let $d\geq 1$ be a square-free integer, set 
$\omega=i\sqrt{d}$ if $-d\equiv 2,3\; (\bmod{4})$, and 
$\omega=(1+i\sqrt{d})/2$ if $-d\equiv 1\; (\bmod{4})$.

Consider the order
$\O_{d}=\Z[1,\vartheta,\vartheta^2,\omega,\omega\vartheta,\omega\vartheta^2]$
of the field $K_d=\Q(\vartheta, i\sqrt{d})$.
If $4d$ (respectively $d$) is co-prime to 49, then $\O_{d}$ is the ring of integers
of $K_d=\Q(\vartheta,i\sqrt{d})$.
We have 

\begin{theorem}
\label{th3}
The order $\O_{d}$ is only monogenic for $d=1$, in which case all generators 
of power integral bases are of the form
\[
\alpha=x_0+x_1\vartheta+x_2\vartheta^2+iy_0+iy_1\vartheta+iy_2\vartheta^2,
\]
where $x_0\in\Z$,
\[
(x_1,x_2,y_0,y_1,y_2)=
(0, 0, 1, -1, 0),
(0, 0, 1, 0, -1),
(0, 0, 2, 0, -1),
\]
\[
(0, 0, 1, 1, -1),
(0, 0, 2, 1, -1),
(0, 0, 0, 1, 0).
\]
\end{theorem}

\noindent
{\bf Proof of Theorem \ref{th3}}.\\

I.A. Assume $-d\equiv 2,3\; (\bmod{4})$ and $d>1$. 
Then we have $y_1=y_2=0$ from (\ref{y12}) and $y_0=\pm 1$ from (\ref{L}).
To calculate $x_1,x_2$ we solve equation (\ref{x12}), that is
\[
|N_{L/Q}(x_1-(-1+\vartheta)x_2)|\leq 1
\]
by Magma. Substituting all possible $x_1,x_2,y_0,y_1,y_2$ into $F(x_1,x_2,y_0,y_1,y_2)$,
we obtain cubic polynomials in $d$ with all negative coefficients,  that is 
(\ref{F3}) can not be satisfied.

I.B. Let $d=-1$. We calculate all solutions $y_1,y_2$ of (\ref{y12}).
For all these explicit values we determine $y_0$ from (\ref{L}).
We determine the possible values of $x_1,x_2$ from (\ref{x12}).
Calculating the indices of $\alpha$ of (\ref{alpha}) for all possible
$x_1,x_2,y_0,y_1,y_2$ we obtain the generators of power integral bases.

II.A. Assume $-d\equiv 1\; (\bmod{4})$ and $d>3$.
Then we have $y_1=y_2=0$ from (\ref{y128}) and $y_0=\pm 1$ from (\ref{L}).
To calculate $x_1,x_2$ we solve equation (\ref{x128}).
We substitute all possible $x_1,x_2,y_0,y_1,y_2$ into $F_3(x_1,x_2,y_0,y_1,y_2)$
and obtain that (\ref{F3}) can not be satisfied.

II.B. Let $d=3$. Using Magma we calculate the possible solutions of (\ref{y128}).
For all these $y_1,y_2$ we calculate $y_0$ from (\ref{L}).
Further, for all $y_1,y_2$ we calculate the solutions $x_1,x_2$ of (\ref{x128}).
Testing the indices of $\alpha$ of (\ref{alpha}) for all these $x_1,x_2,y_0,y_1,y_2$,
we do not get any elements of index 1.

\section{Computational aspects}

All calculations connected with the above examples 
were performed in Maple \cite{maple},
except for solving the Thue equations, which was done
in Kash \cite{kant} and Magma \cite{magma}.
Our procedures were executed 
on an average laptop running under Windows. The CPU time took 
all together some seconds.

\end{document}